\documentclass[reqno]{amsart}

\usepackage{amssymb,amscd,enumerate,graphicx}

\theoremstyle{plain}
\newtheorem{theorem}{Theorem}
\newtheorem{lemma}{Lemma}
\newtheorem{proposition}{Proposition}
\newtheorem{corollary}{Corollary}

\theoremstyle{definition}
\newtheorem{definition}{Definition}

\theoremstyle{remark}
\newtheorem*{remark}{Remark}
\newtheorem*{note}{Note}

\newcommand{\sfi}{\hbox{\ensuremath{\mathcal{S}(\varphi)}}}

\newcommand{\N}{\hbox{\ensuremath{\mathbb{N}}}}
\newcommand{\Z}{\hbox{\ensuremath{\mathbb{Z}}}}

\newcommand{\R}{\hbox{\ensuremath{\mathbb{R}}}}
\newcommand{\C}{\hbox{\ensuremath{\mathbb{C}}}}
\newcommand{\ggr}{\hbox{\ensuremath{\mathbf{g}}}}
\newcommand{\hh}{\hbox{\ensuremath{\mathbf{h}}}}
\newcommand{\tk}{t_0}
\newcommand{\BB}{\hbox{\ensuremath{\mathcal{B}}}}

\newcommand{\comment}[1]{}

\DeclareMathOperator{\kker}{Ker}
\DeclareMathOperator{\supp}{Supp}
\DeclareMathOperator{\sspan}{span}

\numberwithin{theorem}{section}
\numberwithin{proposition}{section}
\numberwithin{lemma}{section}
\numberwithin{corollary}{section}
\numberwithin{definition}{section}
\numberwithin{equation}{section}

\begin{document}



\title{Refinable Shift Invariant Spaces  in $\R^d$}

\author[C.A.Cabrelli]
{Carlos A.~Cabrelli}
\address{
Depto. de
Matem\'atica \\ FCEyN\\ Univ.
de Buenos Aires\\ Cdad. Univ., Pab. I\\ 1428 Capital
Federal\\ ARGENTINA\\ and CONICET, Argentina}
\email[%
]{cabrelli@dm.uba.ar, sheinek@dm.uba.ar, umolter@dm.uba.ar}

\author[S.B.Heineken]{Sigrid B.~Heineken}

\author[U.M.Molter]{Ursula~M.~Molter}

\thanks{The research of
the authors is partially supported by Grants:
CONICET PIP456/98, and UBACyT X610}



\date{\today}


\begin{abstract}
Let $\varphi: \R^d \longrightarrow \C$ be a compactly supported
function which satisfies a refinement equation of the form
\begin{equation*}
\varphi(x) = \sum_{k\in\Lambda} c_k \varphi(Ax - k),\quad
c_k\in\C,
\end{equation*}
where $\Gamma\subset\R^d$ is a lattice, $\Lambda$ is a finite
subset of $\Gamma$, and $A$ is a dilation matrix. We prove, under
the hypothesis of linear independence of the $\Gamma$-translates
of $\varphi$, that there exists a correspondence between
the vectors of the Jordan basis of a
finite submatrix of $L=[c_{Ai-j}]_{i,j\in\Gamma}$
and a finite dimensional subspace
$\mathcal H$ in the shift invariant space
generated by $\varphi$.
We provide a basis of $\mathcal H$ and show that its elements
satisfy a property of homogeneity associated to the eigenvalues of
$L$.
If the function $\varphi$ has accuracy $\kappa$, this basis can be
chosen to contain
a basis for all the multivariate polynomials of degree less than
$\kappa$.
These latter functions are associated to eigenvalues that are powers
of
the eigenvalues of $A^{-1}$. Further we show that the dimension of
$\mathcal H$
coincides with the local dimension of $\varphi$, and hence, every function in the shift invariant space generated by $\varphi$ can be written locally as a linear combination of translates of the homogeneous functions.
\end{abstract}

\maketitle


\section{Introduction}
Let $\Gamma$ be a lattice, i.e. $\Gamma$ is the image of $\Z^d$
under any nonsingular linear transformation and $A$ be a {\em
dilation matrix} associated to $\Gamma$, (i.e.
$A(\Gamma)\subset\Gamma$ and all eigenvalues of $A$ satisfy
$|\lambda|>1$). We will say that a compactly supported function
$\varphi: \R^d \longrightarrow \C$ is {\em refinable} with respect
to $A$ and $\Gamma$, if it satisfies the {\em dilation equation}
\begin{equation}\label{dil-eq}
\varphi(x) = \sum_{k\in\Lambda} c_k \varphi(Ax - k),\quad x\in
\R^d,
\end{equation}
for some finite subset $\Lambda\subset\Gamma$, and coefficients
$c_k\in\C$.

The {\em Shift Invariant Space} (SIS) generated by $\varphi$ is
the space
$$ {\mathcal S}(\varphi) = \left\{f: \R^d \longrightarrow \C: f(x)
= \sum_{k\in \Gamma} y_k \varphi(x+k), \quad
{y_k \in \C, k\in\Gamma}\right\}. $$
Note that since $\varphi$ is compactly supported, the right hand side of the
previous equation is well defined.
Even though the SIS is an infinite dimensional space, the fact that
its
generator is compactly supported yields ``locally'' a finite number
of generators. More precisely,  let $E \subset \R^d$ be a fundamental
domain for the lattice $\Gamma$, and
\begin{equation}\label{efi}
E(\varphi) = \{f/_E: E \longrightarrow \C: f/_E(x) = f(x)\ \forall x
\in E,
f \in \sfi \}.
\end{equation}
The space $E(\varphi)$ is finite dimensional. A canonical
set of generators is the set
$$
 \{\varphi(x-k)/_E: E \longrightarrow \C;\,k  \text{ such that }
|(\supp(\varphi) + k)\cap E| > 0 \}.
$$
The algebraic dimension of the vector space $E(\varphi)$ will be
called the
{\em local
dimension} of $S(\varphi)$ and a basis of $E(\varphi)$
a local basis for  $S(\varphi)$.

If  $\varphi$ satisfies a refinement equation~\eqref{dil-eq},
one may not explicitly
know the function, but from properties of the coefficients of
~\eqref{dil-eq} one often can deduce (properties) of these finite
generators.
One  question is if it is possible to obtain a different set of
generators
with specific properties. In particular it is important to know
 if some (or any) of these generators can be chosen
to be polynomials. The {\em accuracy} of $\varphi$ is the maximum
integer
$\kappa$, such that all polynomials of degree less or equal to
$\kappa-1$
are contained in $S(\varphi)$. Hence, if $\varphi$ has accuracy
$\kappa$,
one can choose a local basis containing
$\alpha_{\kappa}=\sum_{s=0}^{\kappa-1} d_s$
linearly independent polynomials, where $d_s$ is the number of
linearly
independent monomials of degree $s$.

In the one dimensional case, with dilation $2$, the generator
$\varphi$  satisfies $\varphi(x) = \sum_{k=0}^N c_k \varphi(2x-k)$.
In this case  $\alpha_{\kappa} = \kappa$ and the accuracy
of $\varphi$ is related to spectral properties of a finite matrix $T$.

Precisely, under the hypothesis of linear independence
of the integer translates of the generator,
 $\varphi$ has accuracy $\kappa$, if and only if
$\{1,2^{-1},...,2^{-(\kappa-1))}\}$ are
 eigenvalues of the $(N+1)\times (N+1)$ matrix $T$ defined by
$T = \{c_{2i-j}\}_{i,j=0,...,N}$, (the {\em scale matrix}) {\bf
and}
there exist polynomials $p_0,...,p_{\kappa-1}$ of degree$(p_i)=i$
such that each of the the vectors $v_i= \{p_i(k)\}_{k=0,...,N}$
is a left eigenvector of $T$ corresponding to the eigenvalue
$2^{-i}$ (\cite{Dau88}, \cite{CHM98}).
Here the fact that the powers of $1/2$ are eigenvalues
of T is related to the dilation factor $2$ in the dilation
equation.

So, if $\varphi$ has accuracy $\kappa$, we know $\kappa$ linearly
independent functions in $E(\varphi)$. This set of functions
can then be extended to a local basis of $S(\varphi)$.
A natural question is, in which way can this completion be done?
If some eigenvalues are associated to nice local basis functions,
(i.e. local polynomials) would it be possible to extend this set of
nice
functions to a basis of  $E(\varphi)$ using the remaining eigenvalues?
What properties will these new functions have?

Blu and Unser~\cite{BU02} in the study of radial basis functions,
and later Zhou~\cite{Zho02} gave the first clue for the answer to
these questions.
They showed that associated to an arbitrary eigenvalue $\lambda$
of the matrix $T$ there is a function in $S(\varphi)$ that satisfies
that
$h(2x) = \frac{1}{\lambda} h(x)$. In particular, the monomial $x^k$,
satisfies this property, for $\lambda = 2^{-k}$.
So the set of functions associated in this way to all the eigenvalues
of $T$ is a linearly independent set.
Then in the case that the matrix $T$ is diagonalizable it is possible
to complete
the local basis for $S(\varphi)$.

In~\cite{CHnM05a} the problem was completely solved. They showed
that to each vector from a basis that gives the Jordan form of the
matrix $T$ it is possible to associate a function $h$ in
$S(\varphi)$ that they called $(\lambda,r)$-homogeneous  and
satisfies that $$ (D_2-\lambda I)^rh=0. $$ Here $D_2f$ is the
dilation operator defined by $D_2f(x)= f(\frac{x}{2})$, and
$\lambda$ is an eigenvalue of $T$. In particular, the functions
associated to eigenvectors are $(\lambda,1)$-homogeneous, and
correspond to the one's obtained before. These
$(\lambda,r)$-homogeneous functions in $S(\varphi)$
 are linearly independent and provide a
local basis. This local basis contain all the  monomials
$x^k$  within the accuracy.
The generator $\varphi$ can be completely obtained from this local
basis.

The goal of this paper is to carry on this study to $\R^d$, with a
general
dilation matrix and an arbitrary full rank lattice.

When moving to higher dimensions, the situation turns much more
complicated. In analogy to the one dimensional case, we will
consider functions  that satisfy the  relation $h(Ax) =
\frac{1}{\lambda}h(x)$ (see section~\ref{4}). Here $h:\R^d
\rightarrow \C$, $A$ is a $d\times d$ invertible matrix and
$\lambda \in \C$. More in general, we will consider functions
satisfying that $(D_A-\lambda I)^rh=0$. To avoid any ambiguity, we
will say that functions satisfying this equation are in the  class
$\mathcal{H}(A,\lambda ,r)$, in place to use the word homogeneous,
since we will also be dealing with polynomials that are
homogeneous in the standard way, (i.e. a polynomial $p$ of degree
$s$ is homogeneous, if $p(ax)= a^s p(x), x\in \R^d
 \,\forall a \in \R)$. Note however, that with
this definition, for $d=2$ the monomial $h(x_1,x_2) = x_1x_2$ will
be in $\mathcal{H}(A,\lambda,1)$ only if $A$ is diagonal and $\lambda=
\frac{1}{A_{11} A_{22}}$.

When trying to extend the notion of accuracy from one to higher
dimensions, it became apparent that the fact that $\varphi$ has
accuracy $\kappa$ is not immediately related to spectral
properties of a finite submatrix of $L = [c_{Ai - j}]_{i,j\in
\Gamma}$. The relation is much more subtle and involved (see
\cite{CHM98}, \cite{CHM99}, \cite{CHM00}).

In spite of this, however, in this paper we are able to obtain an
analogous result to the 1-dimensional case. Again we are able to
show that a local basis of ${\mathcal S}(\varphi)$ can be obtained
using solely functions from
 $\mathcal{H}(A,\lambda,r)$, where $\lambda$ is an eigenvalue for
a finite submatrix $T$ of $L$. The result is very pleasing, in the
sense that the 1-dimensional results are completely recovered and
one obtains a {\em different} way of writing the functions of
$\sfi$, namely, each function in $\sfi$ can be written locally as
a linear combination of the translates of functions in
$\mathcal{H}(A,\lambda,r)$ (cf. equation~\eqref{f-homog}).

In particular,  if $\varphi$ has accuracy $\kappa$, then we will
find $\alpha_{\kappa}$ linearly independent polynomials that are
in the class $\mathcal{H}(A,\lambda,r)$ for some eigenvalue
$\lambda$ of $T$, and some $r \in \N $.

The difficulty here is to find the appropriate matrix $T$.
Since we are in $\R^d$, the indexes vary along a $d$-dimensional
lattice, so to write $L$ as a matrix, one has to {\em order}
the points. Which order is not important, as long as it
is always the same. In the one dimensional case, it was
straightforward
to look at a submatrix of $L$ that was intimately related to
the support of $\varphi$. In the higher dimensional setting, it may be a difficult problem
to determine the support exactly.
This is one of the problems one has to overcome to solve the
question raised here.

The paper is organized as follows: In section~\ref{s2} we briefly
review some geometric properties of the support of a refinable
function $\varphi$ related to the dilation $A$ and define the
finite matrix  $T$ whose spectral properties will be fundamental
for our analysis of the class $\mathcal{H}(A,\lambda,r)$. In
section~\ref{3}, we relate the spectral properties of the infinite
matrix $L$ to those of $T$, and in section~\ref{4} we define the
class $\mathcal{H}(A,\lambda,r)$ and show how one associates one
of these functions to each vector of the Jordan basis of  $T$. We
further prove that these functions are a basis for the space
$E(\varphi)$~\eqref{efi}. Finally, in section~\ref{5}, if
$\varphi$ has accuracy $\kappa$, using results from~\cite{CHM03},
we show that the space of all functions in the class
$\mathcal{H}(A,\lambda,r)$ in $S(\varphi)$, contains $\alpha_k$
linearly independent polynomials.

\section{Attractors, Tiles and Admissible Sets}
\label{s2}

Let $\Gamma$ be a lattice and $A$ a dilation matrix associated to
$\Gamma$.
Then $A$ has
integer determinant and the group $\Gamma/A(\Gamma)$ has order
$|\det(A)|$ (see for example \cite{Woj97}).
Set
\begin{equation*}
m=|\det(A)|,
\end{equation*}
and let $D=\{d_1,\ldots,d_m \}$ be a set of representatives of the
group $\Gamma/A(\Gamma)$ of order $m$. We call $D$ a {\em full set
of digits}, or {\em digit set}.

The cosets
\begin{equation*}
\Gamma_i=A(\Gamma)-d_i=\{Ak-d_i:k\in\Gamma\},\quad d_i\in D,
\end{equation*}
form a partition of $\Gamma$. Assume that
$\gamma_1,\ldots,\gamma_d$ is a set of generators for $\Gamma$,
that is, $\gamma_1,\ldots,\gamma_d$ are linearly independent
vectors in $\R^d$ and
\begin{equation*}
\Gamma=\{l_1\gamma_1+\ldots+l_d\gamma_d: l_i\in\Z\}.
\end{equation*}
We will call the set
\begin{equation*}
P=\{x_1\gamma_1+\ldots+x_d\gamma_d: 0\leq x_i<1\}
\end{equation*}
a {\em fundamental domain} for the group $\R^d/\Gamma$.

\subsection{Attractors}
For each $k\in\Gamma$, we define $w_k:\R^d\longrightarrow \R^d$ by
\begin{equation*}
w_k(x)=A^{-1}(x+k).
\end{equation*}
Since $A$ is a dilation matrix, $A^{-1}$ is contractive for some
appropriate norm in $\R^d$, so each $w_k$ is a contractive mapping
on $\R^d$ for that norm.

The space
\begin{equation*}
\mathcal{H}(\R^d)=\{K\subset\R^d: K\neq\emptyset \text{ and $K$ is
compact}\},
\end{equation*}
is a complete metric space under the Hausdorff metric $d$
defined by

$$d(B,C)=\text{inf}\left\{\varepsilon > 0: B\subset
C_{\varepsilon}\text{ and }C\subset B_{\varepsilon}\right\},$$
where $$B_{\varepsilon}=\left\{x\in \R^d :\text{ dist}(x,B)
<\varepsilon\right\}.$$

For each finite subset $H\subset\Gamma$, we define
\begin{equation*}
  w_H(B)=\bigcup_{k\in H}w_k(B)=A^{-1}(B+H)
\end{equation*}
It can be shown that $w_H$ is a contractive map in
$\mathcal{H}(\R^d)$ (using that each $w_k$ is a contractive
mapping on $\R^d$). Consequently, by the Contraction Mapping
Theorem, there exists a unique nonempty compact set
$K_H\subset\R^d$ such that
\begin{equation*}
w_H(K_H)=K_H,\quad i.e.\quad K_H=A^{-1}(K_H+H).
\end{equation*}
In fact, we can write
\begin{equation}\label{attractor}
K_H=\sum_{j=1}^\infty A^{-j}(H)=\left\{\sum_{j=1}^\infty
A^{-j}h_j:h_j\in H\right\}.
\end{equation}
The set $K_H$ is called the {\em attractor} of the iterated
function system generated by $\{w_k\}_{k\in H}$.
\cite{Hut81}.
\subsection{Tiles}
Given the digit set $D=\{d_1,\ldots,d_m \}$, we consider the
attractor
\begin{equation}\label{Q}
Q=K_D=\sum_{j=1}^\infty A^{-j}(D)=\left\{\sum_{j=1}^\infty
A^{-j}\varepsilon_j:\varepsilon_j\in D\right\}
\end{equation}
of the iterated system generated by $\{w_d\}_{d \in D}$. We have
that, for $\gamma\in\Gamma$
\begin{equation}\label{tra}
K_{D+\gamma}=\sum_{j=1}^\infty A^{-j}(D+\gamma)=\sum_{j=1}^\infty
A^{-j}D+(A-I)^{-1}\gamma=K_D+(A-I)^{-1}\gamma.
\end{equation}
Therefore we can assume without loss of generality that $0\in D$,
and hence, by equation (\ref{Q}), we have $0\in Q$.

The set $Q$ satisfies the following  properties (see \cite{Ban91}
and \cite{GM92}):

\begin{itemize}
{\em
\item[a)] $\bigcup_{k\in\Gamma}Q+k=\R^d$.
\item[b)] $Q^0\neq\emptyset,  Q=\overline{Q^0}$, and $|\partial
Q|=0$.
\item[c)] $|Q\cap (Q+k)|=0$ for every $k\in\Gamma-\{0\}$ if and only
if
$|Q|=|P|$, where $P$ is a fundamental domain for $\R^d/\Gamma$. In
this case $Q\cap(Q+k)\subset\partial Q$ for all
$k\in\Gamma-\{0\}$.}
\end{itemize}

A longstanding problem was the question of whether for each
dilation matrix $A$ there exists a full set of digits $D$ such
that the corresponding attractor $Q$ is a tile. A counterexample
was found recently. (See \cite{Pot97} and also \cite{LW99}.) We
will assume in this paper that $Q$ is a {\em tile}, or a {\em
fundamental domain for $\Gamma$}, that is, the $\Gamma$-
translates $\{Q+k\}_{k \in \Gamma}$ cover $\R^d$ with overlaps of
measure zero (hence $|Q|=|P|$). Then the local dimension of
$\sfi$, will be the dimension of $Q(\varphi)$. See  \eqref{efi}.

\subsection {Admissible sets}
\label{admisi}
Let $H$ be a fixed finite subset of $\Gamma$
\begin{definition}
We say that a set $\Omega\subset\Gamma$ is $H$-{\em admissible} if
\begin{equation}\label{adm}
  A^{-1}(\Omega+H)\cap \Gamma\subset\Omega,
\end{equation}
which is equivalent to say that
$w_H(\Omega)\cap\Gamma\subset\Omega$.
\end{definition}
\begin{remark} \label{subset} If $H \subset H^{'}$ and $\Omega$ is
$H^{'}$-admissible, then
$\Omega$ is $H$-admissible.
\end{remark}
We immediately have the following Proposition:
\begin{proposition} \label{admis}
If $\Omega_H$ is defined as
$\Omega_H=K_H\cap\Gamma$, then
$\Omega_H$ is an $H$-admissible set.
\end{proposition}
\begin{proof}
Since $\Omega_H\subset K_H$, we have
\begin{equation*}
w_H(\Omega_H)\cap\Gamma\subset w_H(K_H)\cap\Gamma=\Omega_H,
\end{equation*}
which shows the desired property.
\end{proof}
Let $\ell(\Gamma)$ be the space of all sequences defined in
$\Gamma$, and let $L$ be the infinite matrix
associated to the refinement equation~\eqref{dil-eq}, $L_{ij}=
c_{Ai-j}$, whenever $Ai-j \in \Lambda$, and $L_{ij}=0$ in all other
cases.

In this paper we will mainly be interested in $\Lambda$-admissible
sets. The reason for that is that if $\Omega \subset \Gamma$ is
$\Lambda$-admissible, then
the space $\ell(\Omega)
= \{Y \in \ell(\Gamma): y_k = 0, k \not\in \Omega \}$ is
right invariant under $L$.

We will need to ``extend'' finite vectors to infinite
ones with certain prescribed properties, and such that they
coincide with the finite one if restricted to a
finite subset of the lattice. Therefore, the following
Proposition found in \cite{CHM04}, will be very useful.

\begin{proposition}[CHM04]\label{seq}
For each finite $H \subset \Gamma$,
there exists a strictly increasing sequence
$\{\Omega_n\}_{n\geq 0}$ of $H$-admissible sets whose union is
$\Gamma$, such that $\Omega_0=\Omega_H$ and
\begin{equation}\label{Inc}
w_H(\Omega_{n+1})\cap\Gamma\subset\Omega_n
\end{equation}
for all $n\geq 0$.
\end{proposition}

\begin{proof}
Let $\|\cdot\|$ be any norm in $\R^d$ such that $\|A^{-1}\|<1$ and
fix $\varepsilon >0$, such that $H\subset B(\varepsilon)$, where
$B(\varepsilon)=\{x\in\R^d:\|x\|\leq \varepsilon\}$,
the closed ball with radius $\varepsilon$ centered at the origin.
Now set
\begin{equation}\label{del}
  \delta_0=\frac{\varepsilon}{\|A^{-1}\|^{-1}-1}.
\end{equation}
Choose $\delta > \delta_0$ in such a way that $\Omega_H\subset
B(\delta)$ and set $F_0=B(\delta)$. Since $\delta > \delta_0$, we
have $\|A^{-1}\|(\delta +\varepsilon)< \delta$. Hence,
\begin{equation*}
w_H(F_0)=A^{-1}(B(\delta) +H)\subset A^{-1}(B(\delta
+\varepsilon))\subset B(\|A^{-1}\|(\delta +\varepsilon))\subset
B(\delta)=F_0.
\end{equation*}
We define recursively $F_{j+1}=w_H(F_j)$ for $j\geq 0$. It is easy
to see, by induction, that $F_{j+1}\subset F_j$ for every $j$.
Since $F_0$ is compact, the Contraction Mapping Theorem tells us
that $\bigcap F_j=K_H$. It follows that $F_j\cap\Gamma=\Omega_H$
for every $j$ large enough, and consequently $\{F_j\cap\Gamma\}$
is a finite collection of sets. Let
\begin{equation*}
\Omega_H=\Omega_0\subsetneq\Omega_1\subsetneq\cdots\subsetneq\Omega_N=F_0\cap\Gamma
\end{equation*}
be the distinct elements of this collection and fix $0\leq n<N$.
Since there exists a $j\in\N$ such that
\begin{equation*}
\Omega_n=F_j\cap\Gamma\subsetneq F_{j-1}\cap\Gamma=\Omega_{n+1},
\end{equation*}
we have
\begin{equation}
w_H(\Omega_{n+1})\cap\Gamma\subset
w_H(F_{j-1})\cap\Gamma=F_j\cap\Gamma=\Omega_n.
\end{equation}
So  inclusion (\ref{Inc}) holds for $n=0,1,\ldots,N-1$.

Now we set $\delta_N=\delta$ and define recursively,
$\delta_{n+1}=\frac{\delta_n}{\|A^{-1}\|}-\varepsilon$ for $n\geq
N$. The sequence of numbers $\delta_N<\delta_{N+1}<\cdots$ is
increasing. Define $\Omega_n=B(\delta_n)\cap\Gamma$ for $n>N$.
If $\Omega_{n+1} = \Omega_{n}$, we skip that one and continue until
$\Omega_{n+k} \not= \Omega_n$. In this way we obtain a strictly
increasing
sequence of sets $\{\Omega_k\}_{k\geq N}$.
Combining with the sets $\Omega_0,\ldots,\Omega_N$ constructed
previously, we have a strictly increasing sequence
$\{\Omega_n\}_{n\geq 0}$. The inclusion
\begin{equation}\label{Inc2}
w_H(\Omega_{n+1})\cap\Gamma\subset\Omega_n
\end{equation}
holds for every $n\in\N_0$, since for $n\geq N$,
again there exist a $j\in\N$ such that
\begin{equation*}
\Omega_n=B(\delta_j)\cap\Gamma\subsetneq
B(\delta_{j+1})\cap\Gamma=\Omega_{n+1},
\end{equation*}
and then
\begin{equation}
w_H(\Omega_{n+1})=A^{-1}(\Omega_{n+1}+H)\subset
B(\|A^{-1}\|(\delta_{j+1}+\varepsilon))=B(\delta_j).
\end{equation}
We already showed that $\Omega_0=\Omega_H$ is $H$-admissible.
Since $\Omega_n\subset\Omega_{n+1}$, it follows from (\ref{Inc2})
that $\Omega_{n+1}$ is $H$-admissible for every $n\in\N_0$, which
completes the proof.
\end{proof}
\begin{corollary}\label{seq-1}
If $H \subsetneq H^{\prime} \subset \Gamma$, then there exists $n_0
\geq 1$ and a strictly increasing
sequence
$\{\Omega_n\}_{n\geq 0}$ of $H$-admissible sets whose union is
$\Gamma$, such that
\begin{equation*}
\Omega_0=\Omega_H, \quad \Omega_{n_0} = \Omega_{H^{\prime}}, \quad
\text{and} \quad w_H(\Omega_{n+1})\cap\Gamma\subset\Omega_n \
\text{for all $n\geq 0$.}
\end{equation*}
\end{corollary}
\begin{proof}
First construct the sequence $\{\Omega^{\prime}_{n}\}$ associated to
$H^{\prime}$ using the previous proposition.
Note that by Remark~\ref{subset} the sets $\Omega^{\prime}_{n}$ are
also $H$-admissible.

Using the notation of the previous proof, let $j_0$ be such that
$F_{j_0} \cap \Gamma = \Omega_{H^{\prime}}$.
Now consider the sequence $\{G_j\}_{j\geq j_0}$, where $G_{j_0} =
F_{j_0}$ and $G_{j+1} = w_H(G_j)$, for
$j\geq j_0$. Since
\begin{equation*}
G_{j_0+1} = w_H(G_{j_0}) = A^{-1}(F_{j_0}+H) \subseteq
A^{-1}(F_{j_0}+H^{\prime}) = w_{H^{\prime}}(F_{j_0}) \subseteq
F_{j_0} = G_{j_0},
\end{equation*}
then $G_{j+1} \subseteq G_{j}$ and therefore $\cap_{j\geq j_0} G_j =
K_H$ and hence,
$\{G_j \cap \Gamma\}$ is again a finite collection of sets of say
$n_0+1$ elements. Consider now the distinct
elements
\begin{equation*}
\Omega_H = \Omega_0 \subsetneq \Omega_1 \subsetneq \dots \subsetneq
\Omega_{n_0} = F_{j_0} \cap \Gamma =
\Omega_{H^{\prime}},
\end{equation*}
and let
$$ \Omega_{n_0 +k} = \Omega^{\prime}_k.$$
This new sequence satisfies all the desired properties.
\end{proof}

For a more complete treatment of admissible sets see \cite{CHM04}
and also Jia~\cite{Jia98}.

\section{Spectral Properties of $L$}
\label{3}

Let us now return to the refinement equation~\eqref{dil-eq},
$\varphi(x) = \sum_{k \in \Lambda} c_k \varphi(Ax-k)$.
If  we consider the infinite column vector
\begin{equation}
\Phi(x)=\{\varphi(x+k)\}_{k\in\Gamma},
\end{equation}
this equation becomes
\begin{equation}
\Phi(x)=L\Phi(Ax).
\end{equation}
It can be shown (see \cite{CHM00}) that the set
$K_\Lambda$, which is the particular case taking  $H=\Lambda$
in \eqref{attractor}, satisfies that if $\varphi$ is a
compactly supported solution of the refinement equation
(\ref{dil-eq}), then $\supp(\varphi)\subset K_\Lambda$.
Also, by Proposition~\ref{admis}, the set
$\Omega_{\Lambda}=K_{\Lambda}\cap\Gamma$
is $\Lambda$-admissible. However, it is not necessarily true that
$\supp \varphi \subset \cup_{\lambda \in \Omega_{\Lambda}} Q
+ \lambda$.

We will therefore consider the bigger set $\Omega^{\prime}=
K_{\Lambda^{\prime}}\cap \Gamma$, where $\Lambda^{\prime} = \Lambda -
D \supset \Lambda$. In~\cite{CHM04} it was shown that the
translations of $Q$ using all elements of $\Omega^{\prime}$
cover the support of the compactly solution to~\eqref{dil-eq}.
Moreover, $\Omega^{\prime}$ is $\Lambda^{\prime}$
admissible, and hence
also $\Lambda$-admissible.
As noted earlier, the $\Lambda$-admissibility of
$\Omega^{\prime}$ guarantees that the space
$\ell(\Omega^{\prime})
= \{Y \in \ell(\Gamma): y_k = 0, k \not\in \Omega^{\prime}
\}$ is
right invariant under $L$.

Let now $\{\Omega_n\}_{n\geq 0}$ be a sequence of subsets of $\Gamma$ that satisfies:
\begin{itemize}
\item $\Omega_0 = \Omega_{\Lambda}$
\item For $i \geq 0$, $\Omega_i \subsetneq \Omega_{i+1},$
 and $\cup_i \Omega_i = \Gamma$
\item For $i \geq 0$ $\Omega_i$ are $\Lambda$-admissible and
$w_{\Lambda}(\Omega_{i+1}) \cap \Gamma \subseteq \Omega_i.$
\item $\Omega_{n_0} = \Omega^{\prime}$
\item For $i \geq n_0$, $\Omega_i$ are $\Lambda^{\prime}$-admissible and
$w_{\Lambda^{\prime}}(\Omega_{i+1}) \cap \Gamma \subseteq \Omega_i.$
\end{itemize}
These sets exist by Proposition~\ref{seq} and its Corollary~\ref{seq-1}.

We denote by $\{T_n\}_{n\geq 0}$ the finite
submatrices of $L$
\begin{equation}
T_n =[c_{Ai-j}]_{i,j\in\Omega_{n}}.
\end{equation}
Since $\Omega_{n} \subset \Omega_{n+1}$, if the order in $\Gamma$ is appropriately chosen,
actually $T_n$ is a submatrix of $T_{n+1}$, for each $n$.

Let $Y=\{y_k\}_{k\in\Gamma}\in\ell(\Gamma)$ be an infinite row
vector, and
$P_n:\ell(\Gamma)\longrightarrow\C^{1\times\Omega_{n}}, n\geq 0$
be the
restriction mappings defined by
\begin{equation} \label{projection}
P_nY=\{y_k\}_{k\in\Omega_{n}}.
\end{equation}
We consider $L-\lambda I:\ell(\Gamma)\longrightarrow\ell(\Gamma)$
the left-multiplication operator who maps $Y\longrightarrow
Y(L-\lambda I)$ (where $I$ is the identity operator acting on
$\ell(\Gamma))$. By abuse of notation, $I$ will be any identity
operator, no matter on which space it is acting on.
\begin{note}
In what follows we will use powers of the matrix $(L-\lambda I)$.
Note that these powers are point-wise well defined, since the rows
of the matrix $L$ have a finite number of non-zero elements.
\end{note}

The next proposition shows the relation between the
spectrum of $L$ and $T_n$:
\bigskip
\begin{proposition}\label{spect}
Consider $\lambda\in\C$, $r\in\N$ and $n \geq 0$.
\begin{enumerate}
\item Let $Y\in\ell(\Gamma)$. We have
\begin{equation}
Y \in \kker(L - \lambda I)^r \quad \text{implies} \quad
P_nY\in\kker(T_n - \lambda I)^r.
\end{equation}
Conversely,
\item If $v\in \kker(T_n - \lambda I)^r$ and $\lambda\neq 0$, then we
can extend v  to an infinite row vector $Y_v$ (i.e.
$Y_v\in\ell(\Gamma)$ and $P_nY_v=v$), so that $Y_v\in
\kker(L-\lambda I)^r$.\label{2}
\item If $\lambda\neq 0, Y\neq 0$ and $Y\in \kker(L - \lambda I)^r$,
then $P_nY\neq 0$. In particular the extension in \eqref{2} of $v$ to $Y_v$ is unique.
\end{enumerate}
\end{proposition}

\begin{proof}
\

\begin{enumerate}
\item First note that $j\in \Omega_{n}$ and
$Ai-j\in\Lambda$ implies $i\in\Omega_{n}$. For, in this
case, $Ai\in\Omega_{n}+\Lambda$ and since
$\Omega_{n}$ is a $\Lambda$-admissible set it follows
that $i\in
A^{-1}(\Omega_{n}+\Lambda)\cap\Gamma\subset
\Omega_{n}$. Hence
\begin{equation}\label{***}
j \in \Omega_{n},\ i \not\in \Omega_{n}
\ \Longrightarrow
[L-\lambda I]_{ij} = 0.
\end{equation}
Moreover, we will show by induction on $r$ that,
\begin{equation} \label{induction}
 \text{if} \;
j\in\Omega_{n}\; \text{and}  \; i \not \in \Omega_{n}, \text{then} \;
{[(L-\lambda I)^r]}_{ij}=0.
\end{equation}
\roman{enumii}
\begin{enumerate}
\item   The case $r=1$ is simply \eqref{***}, since we assume that
$c_k=0$ if $k\not\in\Lambda$.

\item Suppose now that \eqref{induction} holds for some fixed $r \geq 1$.
 Using (a), for $j \in \Omega_n$ we have
\begin{align*}
{[(L-\lambda I)^{r+1}]}_{ij}&=\sum_{k\in\Gamma}{[(L-\lambda
I)^r]}_{ik}[L-\lambda I]_{kj}\\
&=\sum_{k\in\Omega_{n}}{[(L-\lambda I)^r]}_{ik}[L-\lambda I]_{kj}.
\end{align*}
Now, if $i\not\in\Omega_{n}$, the inductive hypothesis
yields
that the last sum is zero.
\end{enumerate}
\noindent{Therefore}, the statement is true for all $r\in\N$.

To prove the first part of the Proposition, let $Y\in\ell(\Gamma)$ and
$Y\in\kker(L-\lambda I)^r$.

Applying the preceding equality, we obtain for each $j \in
\Omega_{n}$
\begin{align*}
\left[(P_nY)(T_n - \lambda I)^r\right]_j &=
\sum_{i\in\Omega_{n}}y_i{[(L-\lambda
I)^r]}_{ij}\\
&=\sum_{i\in\Gamma}y_i{[(L-\lambda I)^r]}_{ij} \\
&= \left[Y(L-\lambda I)^r\right]_j =0.
\end{align*}
This completes the proof of (1).

\item Assume that $v\in\C^{1\times\Omega_{n}},
\lambda\neq 0$ and $v\in\kker(T_n - \lambda
I)^r$. We want to construct a vector $Y_v\in\ell(\Gamma)$ such
that $P_nY_v=v$ and $Y_v\in\kker(L-\lambda I)^r$.

We now prove by
induction on $r$ that,
\begin{equation}
\text{for}\; j \in \Omega_{k+1}, \;  i\not \in\Omega_k,  k \geq 0,  \ {[(L-\lambda I)^r]}_{ij}=
\begin{cases}
0 & \text{for } i\neq j\\ (-\lambda)^r & \text{for } i=j.
\end{cases}\label{ind2}
\end{equation}

\begin{enumerate}
\item Case  $r=1$.
If $i$ were such that $Ai-j \in \Lambda$, then $i \in A^{-1}(\Omega_{k+1}+\Lambda)
\cap\Gamma\subset\Omega_k$. Hence, if $i\not\in\Omega_k$ then $Ai-j\not\in\Lambda$ and so
$L_{ij}=0$, and therefore
$[L-\lambda I]_{ij}=0$ for $i\neq j$, and $[(L-\lambda
I)]_{jj}=-\lambda$.

\item Assume that (\ref{ind2})
holds for $r\geq 1$. Then for  $j \in \Omega_{k+1}$ and
$i\not \in\Omega_k,  k \geq 0,$  we have
\begin{align*}
{[(L-\lambda I)^{r+1}]}_{ij}&=\sum_{\ell\in\Gamma}{[(L-\lambda
I)^r]}_{i\ell}[L-\lambda I]_{\ell j}\\
&=\sum_{\ell\not\in\Omega_k}{[(L-\lambda I)^r]}_{i\ell}[L-\lambda I]_{\ell j},
\end{align*}
since by the inductive hypothesis ${[(L-\lambda I)^r]}_{i\ell}=0$ if
$\ell\in\Omega_k\subset\Omega_{k+1}$ and $i\not\in\Omega_k$. It
follows using the case $r=1$ that
\begin{equation}
{[(L-\lambda I)^{r+1}]}_{ij}=
\begin{cases}
0 & \text{for } i\neq j\\ (-\lambda)^{r+1} & \text{for } i=j,
\end{cases}
\end{equation}
which proves (\ref{ind2}), for every $r \in \N$.
\end{enumerate}
Define now $y_j=v_j$ for $j\in\Omega_{n}$, and define
recursively, for $j\not\in\Omega_{n}$,
\begin{equation}\label{defyv}
y_j=\frac{-1}{(-\lambda)^r}\sum_{i\in\Omega_k}y_i{[(L-\lambda
I)^r]}_{ij},\quad j\in\Omega_{k+1}\backslash\Omega_k, k\geq n.
\end{equation}
The vector $Y_v=\{y_j\}_{j\in\Gamma}$ is an extension of $v$. To see that
 $Y_v\in\kker(L-\lambda I)^r$, since $\left(Y(L-\lambda I)^r\right)_j =
\sum_{i\in \Gamma}
y_i [(L -\lambda I)^r]_{ij}$, we have:
\begin{itemize}
\item If
$j\in\Omega_{n}$, then by \eqref{***}
\begin{align}
\sum_{i\in\Gamma}y_i{[(L-\lambda I)^r]}_{ij}
& =\sum_{i\in\Omega_{n}}y_i{[(T_n - \lambda I)^r]}_{ij} \notag\\
& =\sum_{i\in\Omega_{n}}v_i{[(T_n - \lambda I)^r]}_{ij}=0.
\end{align}
\item If $j\not\in\Omega_{n}$, then there exists
$k\in\N_0, k \geq n$ such
that $j\in\Omega_{k+1}\setminus\Omega_k$. Therefore,
\begin{align*}
\sum_{i\in\Gamma}y_i{[(L-\lambda
I)^r]}_{ij}&=\sum_{i\in\Omega_k}y_i{[(L-\lambda
I)^r]}_{ij}+\sum_{i\not\in\Omega_k}y_i{[(L-\lambda I)^r]}_{ij}\\
&=\sum_{i\in\Omega_k}y_i{[(L-\lambda I)^r]}_{ij}+y_j(-\lambda)^r) \\
& =0.
\ \text{(by~\eqref{defyv})}
\end{align*}
\end{itemize}

\item For the last part of the Proposition, assume that $\lambda\neq 0,
Y\neq 0$, and $Y\in\kker(L-\lambda I)^r$. To show that $P_nY\neq
0$, take $k_0\in\Gamma$ such that $y_{k_0}\neq 0$. If
$k_0\in\Omega_{n}$, there is nothing to prove. Otherwise, let
\begin{equation}
t_0=\min\{k\in \N:k_0\in\Omega_k\}.
\end{equation}
Since $Y(L-\lambda I)^r=0$,
\begin{align*}
\sum_{i\in\Gamma}y_i{[(L-\lambda
I)^r]}_{ik_0}&=\sum_{i\in\Omega_{\tk-1}}y_i{[(L-\lambda
I)^r]}_{ik_0} +\sum_{i\not\in\Omega_{\tk-1}}y_i{[(L-\lambda
I)^r]}_{ik_0}\\&=\sum_{i\in\Omega_{\tk-1}}y_i{[(L-\lambda
I)^r]}_{ik_0}+y_{k_0}(-\lambda)^r=0.
\end{align*}
So, there exist $k_1\in\Omega_k, 0<k<\tk$, such that $y_{k_1}\neq
0$. If $k_1\in\Omega_{n}$, we can stop here. If not, we repeat the
procedure until $k_j\in\Omega_{n}$.
\end{enumerate}
\end{proof}
\begin{remark}
\

\begin{itemize}
\item Since the previous Proposition is true  for {\bf any} set of the sequence
$\Omega_n$, in fact the smallest matrix $T_0$ already has all the
spectral information of $L$.
\item The extension of the
vectors of $\kker(T_0 -\lambda I)^r$ to vectors of $\kker(L-\lambda
I)^r$ will produce {\em intermediate} vectors of $\kker(T_n - \lambda
I)^r$, by the construction of the sets $\Omega_n$ produced in
Corollary~\ref{seq-1}.
\end{itemize}
\end{remark}
For the special case $\lambda=0$, under some
mild assumptions, we have an additional property.
We say that the $\Gamma$ translates
$\{\varphi(\cdot-k)\}_{k\in\Gamma}$ are {\em linearly
independent}, if for any sequence $\{\alpha_k\}_{k\in\Gamma}$ in
$\ell(\Gamma)$,
\begin{equation*}
\sum_{k\in\Gamma}\alpha_k\varphi(\cdot-k)\equiv 0 \quad\text {
implies }\quad\alpha_k=0.
\end{equation*}

\begin{lemma}\label{lemma}
If $\{\varphi(\cdot-k)\}_{k\in\Gamma}$ are linearly independent,
then the operator $L:\ell(\Gamma)\longrightarrow\ell(\Gamma),
Y\longmapsto YL$ is one to one.
\end{lemma}
\begin{proof}
Let $YL=0$. Then
\begin{equation}
Y\Phi(x)=YL\Phi(Ax)=0.
\end{equation}
Since $\{\varphi(\cdot-k)\}_{k\in\Gamma}$ are li\-near\-ly
in\-de\-pen\-dent, $Y\Phi=0$ implies $Y=0$, so $\kker(L)=\nolinebreak
\{0\}$.
\end{proof}

\section{The class $\mathcal{H}(A,\lambda,r)$}
\label{4}

Assume $Y\in\kker(L-\lambda I)^r$. If we define $h(x)=Y\Phi(x)$,
we have
\begin{align*}
0 &= Y (L - \lambda I)^r\Phi(x) = Y \left(\sum_{k=0}^r
\begin{pmatrix}r\\k
\end{pmatrix} (-\lambda)^{r-k} L^{k}\right) \Phi(x) \\
 &= Y \left(\sum_{k=0}^r \begin{pmatrix}r\\k\end{pmatrix}
 (-\lambda)^{r-k} \Phi(A^{-k}x)\right), \\
 &= \sum_{k=0}^r \begin{pmatrix}r\\k\end{pmatrix}
 (-\lambda)^{k} h(A^{k-r}x).
\end{align*}
This leads to the following definition:
\begin{definition}
A function $h$ is in  the class $\mathcal{H}(A,\lambda,r)$, if it
satisfies
\begin{equation}\label{hom}
 \sum_{k=0}^r \begin{pmatrix}r\\k\end{pmatrix}
 (-\lambda)^{k} h(A^{k-r}x)=0 \quad \text{   for every } x\in\R^d.
\end{equation}
\end{definition}

If we define the operator $D_A$ by $D_A(f)(x)=f(A^{-1}x)$, then $h$ is
in $\mathcal{H}(A,\lambda,r)$ if and only if
$$(D_A-\lambda I)^rh=0.$$
A function in $\mathcal{H}(A,\lambda,r)$ will also be said to be of
class
$\mathcal{H}(A,\lambda,r)$.

Note that if $h\in\mathcal{H}(A,\lambda,r)$, then
$h\in\mathcal{H}(A,\lambda,s)$ for every $s\geq r$.
\begin{proposition}
Let $V\subset\R^d$ be a bounded set such that $0\in V^{0}$ and
$V\subset AV$. Set $C = AV\setminus V$ and $a \in  \Z$.
Let $h$ be a function of class $\mathcal{H}(A,\lambda,r)$.
Then the values of $h$ in $\R^d \setminus \{0\}$ can be
determined from its values in any set of the type:
$$
\tilde{C} = \bigcup_{k=a}^{a+r} A^k C.
$$
Furthermore, if $\lambda \not= 1$ then $h(0) = 0$.
\end{proposition}
\begin{proof}
Since $h \in \mathcal{H}(A,\lambda,r)$ we get that
\begin{align}
h(x) &= - \sum_{k=1}^{r}
\left(\begin{array}{l}r\\k\end{array}\right) (-\lambda)^{k}
h(A^kx) \quad \text{and} \label{1}\\
h(x) &= - \sum_{k=1}^{r}
\left(\begin{array}{l}r\\k\end{array}\right) (-\lambda)^{-k}
h(A^{-k}x).\label{e2}
\end{align}
On the other side, it has been proved in \cite{ACM04b} that the
set $C$ satisfies:
\begin{itemize}
{\em
\item[a)]$\bigcup_{j\in\Z}A^{j}C=\R^d\backslash\{0\}$
\item[b)]The sets $\{A^jC\}_{j \in \Z}$ are pairwise disjoint.}
\end{itemize}
So, from \eqref{1} we deduce that the values of $h$ in
$A^{a-1}C$ can be obtained from the values in $\tilde C$, and
analogously, from \eqref{e2}
the values of $h$ in
$A^{a+r+1}C$ can be obtained from the values in $\tilde C$.

Then we proceed inductively to obtain all the values in $\R^d
\setminus
\{0\}$. Finally, it is immediate from the definition, that
$h(0) = 0$ when $\lambda \not= 1$.
\end{proof}

\begin{proposition}\label{li}
Suppose $\{\varphi(\cdot-k)\}_{k\in\Gamma}$ are linearly
independent. Let $f_1, \ldots , f_l \in{\mathcal S}(\varphi),
f_i=Y^i\Phi$, where $Y^i \in \ell(\Gamma)$ . Then $f_1, \ldots ,f_l$ are
linearly independent functions if and only if $Y^1, \ldots ,Y^l$
are linearly independent in $\ell(\Gamma)$.
\end{proposition}
\begin{proof}
Since
\begin{equation}
\sum_{i=1}^l\alpha_if_i=\sum_{i=1}^l\alpha_i(Y^i\Phi)=\left(\sum_{i=1}^l\alpha_iY^i\right)\Phi,
\end{equation}
and the translates of $\varphi$ along the
lattice $\Gamma$ are linearly independent, we conclude that $\sum_{i=1}^l\alpha_if_i\equiv 0$ if
and only if $(\sum_{i=1}^l\alpha_iY^i)=0$, which leads to the desired
result.
\end{proof}
\begin{remark}
Let $E:{\mathcal S}(\varphi)\longrightarrow\ell(\Gamma)$ be the
function that associates to each element of ${\mathcal
S}(\varphi)$, its coordinates in $\{\varphi(x+k)\}$.
Proposition~\ref{li} shows that $E$ is an isomorphism.
\end{remark}
\begin{proposition}\label{corr}
Assume that $\{\varphi(\cdot-k)\}_{k\in\Gamma}$ are linearly
independent.
\begin{enumerate}
\item If $h\in {\mathcal S}(\varphi) , h=Y\Phi$ and $h\in
{\mathcal H}(A,\lambda,r)$, then
$Y \in \kker(L - \lambda I)^r $ and $P_nY \in
\kker(T_n - \lambda I)^r$.

\medskip
\noindent Conversely
\smallskip
\item Assume that $\lambda\neq 0, v\in \kker(T_n - \lambda I)^r$ and
that $Y_v$ is the unique extension of $v$ such that
$Y_v\in\kker(L-\lambda I)^r$ (see Proposition~\ref{spect}). Then
the function $h=Y_v\Phi$ belongs to ${\mathcal H}(A,\lambda,r)$.
\end{enumerate}
\end{proposition}
\begin{proof}
If $h=Y\Phi$ is of class $\mathcal{H}(A,\lambda,r)$, then we have
\begin{align*}
0 & = \sum_{k=0}^r \begin{pmatrix}r\\k\end{pmatrix} (-\lambda)^k
h(A^{k-r} x)\\ & = \sum_{k=0}^r
\begin{pmatrix}r\\k\end{pmatrix} (-\lambda)^kYL^{r-k}\Phi(x)\\ &
= Y \left(L - \lambda I\right)^r \Phi(x).
\end{align*}
Since the $\Gamma$ translates of $\varphi$ are linearly
independent, it follows that $Y(L-\lambda I)^r=0$, and
consequently, by Proposition~\ref{spect}, $P_nY \in
\kker(T_n- \lambda I)^r$.

To prove the second part, note that if $v=0$, the statement is
trivially true. If $\lambda\neq 0$ and
$v\in\kker(T_n- \lambda I)^r$ and $v\neq 0$,
then, by Proposition~\ref{spect} we can extend $v$ to a vector
$Y_v\in\kker(L- \lambda I)^r$, and so the function $h=Y_v\Phi$ is
of class $\mathcal{H}(A,\lambda,r)$.
\end{proof}

\subsection{Jordan decomposition of $T_n$}

 Let $m_n=\#\Omega_{n}$. Consider the set $\Delta_n$ of
eigenvalues of $T_n$ and the associated Jordan
basis $\mathcal{B}_n = \{v_1,\dots,v_{m_n}\}$ of $\C^{m_n}$. For each
$v_i\in\mathcal{B}_n$ we have that
$v_i\in\kker(T_n- \lambda I)^k$ and
$v_i\not\in\kker(T_n- \lambda I)^{k-1}$ for some
$\lambda\in\Delta_n$, and for some $k\geq 1$. So to each $v_i$ there
corresponds a
unique pair $(\lambda,k)$.
Note however that to two different $v_i$s in the basis
there could correspond the same pair.
For each vector of $\mathcal{B}_n$, set $v_i = v_i(\lambda,k)$. If $\lambda\neq
0$, by Proposition~\ref{corr} we can associate to each $v_i(\lambda,k)$
a function $h_{v_i(\lambda,k)}$ in ${\mathcal H}(A,\lambda,k)\cap
{\mathcal S}(\varphi)$.
Since the vectors $v_i(\lambda,k)$ are linearly independent, its
extensions $\{Y_{v_i}\}$ are linearly independent in $\ell(\Gamma)$,
so the functions
$\{h_{v_i(\lambda,k)}\}_{v_i\in\mathcal{B}_n,\lambda\neq 0}$ are
linearly independent.

If $h_1,\ldots,h_l$ are of class $\mathcal{H}(A,\lambda,k_i)$, for
some $k_i, i=1,\dots,l$,
then a
linear combination of them is of class $\mathcal{H}(A,\lambda,k)$,
with
$k=\max_{1\leq i\leq l}(k_i)$,
for if $h=\sum_{i=0}^l\alpha_i h_i(\lambda,k_i)$, then
\begin{equation*}
\left(D_A - \lambda I\right)^k h = \left(D_A - \lambda I\right)^k
\sum_{i=0}^l \alpha_i h_i = \sum_{i=0}^l \alpha_i \left(D_A -
\lambda I\right)^k h_i = 0,
\end{equation*}
and consequently, $h\in{\mathcal H}(A,\lambda,k)$.

Let
\begin{equation} \label{charact}
\chi_{T_n}(x)= \prod_{\lambda \in \Delta_n} (x -
\lambda)^{r_\lambda}
\end{equation}
be the characteristic polynomial of $T_n$, and
set
\begin{equation*}
\mathcal{H}_{\lambda}(\varphi)= \{h \in{\mathcal S}(\varphi): h
\in {\mathcal H}(A,\lambda,k),\ \text{for some $k \geq 1$ }\},
\quad \lambda \in \Delta_n.
\end{equation*}
Then, if $\lambda\neq 0$, using Proposition~\ref{corr},
$\dim(\mathcal{H}_{\lambda}(\varphi))=r_{\lambda}$
and a basis for $\mathcal{H}_{\lambda}(\varphi)$ are the functions of
class
$\mathcal{H}(A,\lambda,k)$ corresponding to the vectors
$v_i\in\mathcal{B}_n$, such that $v_i=v_i(\lambda,k)$, for some $k\geq 1$.
So, if we denote
\begin{equation*}
\mathcal{H} = \bigoplus_{\lambda \in \Delta_n, \lambda \not= 0}
 \mathcal{H}_{\lambda}(\varphi) \subset {\mathcal S}(\varphi),
\end{equation*}
then,  $\dim(\mathcal{H})=m_n-r_0$, where $r_0$ is the dimension of the subspace
generated by  the vectors of the Jordan basis associated to
$\lambda=0$. (Note that $m_n = \sum_{\lambda \in \Delta_n, \lambda \not=0} r_{\lambda}$.)

In order to be able to include the case $\lambda=0$ in our analysis, we need to consider the case in which $\supp \varphi \subset \cup_{\omega \in \Omega_n} Q+\omega$.
This will guarantee, that (except for a possible set of measure zero), $\varphi(x+k) = 0$ if
$k \not\in \Omega_n$.

\subsection{The case in which $\Omega_n$ contains $\supp(\varphi)$} If we recall the choice of the sequence $\{\Omega_n\}$ at the beginning of section~\ref{3}, it is clear, that for $n \geq n_0$, we have that  $\supp \varphi \subseteq \cup_{\omega \in \Omega_n} Q+\omega$, and hence, the local dimension of $\sfi$ is $\dim \sspan\{\varphi(x+k)\}_{k \in \Omega_n}$.

In that case for  $x \in Q^{\circ}$, $\lambda \not=0$, $h_{v_i(\lambda,k)}(x) = \sum_{l \in \Omega_n} [v_i(\lambda,k)]_l\varphi(x+l)$ since $\varphi(x+j) = 0$ if $j \not\in \Omega_n$.

Moreover for $\lambda=0$ we have the following Proposition:
\begin{proposition} \label{lambda0}
Let $n \geq n_0$, and let $r_0$ be the power of $x$ in $\chi_{T_n}$ (c.f.~\eqref{charact}). Consider $v_i = v_i(0,r)$, with $r \leq r_0$.  Define $h(x) =  \sum_{k \in \Omega_n}
[v_i]_k \varphi(x+k)$. Then $h \equiv 0$ a.e. on $Q$.
\end{proposition}

We postpone the proof to remark that with this Proposition, if $n
\geq n_0$, and $\BB_n$ is (as before) the matrix whose rows are
the vectors of the Jordan basis for $T_n$,  then \begin{equation}
\label{homog-fi} \left[\begin{array}{l}h_{v_1(\lambda,k)}(x)\\
\vdots\\ h_{v_{m_n}(\lambda,k)} (x)
\end{array}
\right] = \BB_n P_n\Phi(x) \quad \text{a.e.}\; x\in Q.
\end{equation}
Hence, since the matrix $\BB_n$ is invertible the local dimension
of $\sfi$ coincides with $\dim \sspan\{h_{v_i(\lambda,k)}(x), x
\in Q, v_i \in \BB_n\}$, which is equal to the dimension of
$\mathcal{H}$.

So the local dimension of $\sfi$ can be found by finding the
Jordan form of any of the finite matrices $T_n$ as long as $n \geq
n_0$.

Moreover, any function of the shift-invariant space $\sfi$ can be written as a linear combination of the lattice translates of the homogeneous functions. Namely, let $f\in \sfi$, then
\begin{equation}\label{alpha}
f(x) = \sum_{\gamma \in \Gamma} \alpha_{\gamma} \varphi(x+\gamma) \quad x \in \R^d, \alpha_{\gamma} \in \C.
\end{equation}
If we call $\ggr:\R^n \longrightarrow \C^{m_n}$ and $\hh: \R^n \longrightarrow \C^{m_n}$ the functions
\begin{equation}
\ggr(x) = P_n\Phi(x) \chi_{Q}(x)
\quad \text{and} \quad
\hh(x) = \left[\begin{array}{l}h_{v_1(\lambda,k)}(x) \\
\vdots\\
h_{v_{m_n}(\lambda,k)}(x)
\end{array}
\right] \chi_{Q} (x),
\end{equation}
and for $\gamma \in \Gamma$ we denote by ${\overline{\alpha}}_{\gamma} = (\alpha_{i_1+\gamma}, \dots, \alpha_{i_{m_n} + \gamma})$ the vector of length $m_n$ whose indices are in $\Omega_n+\gamma$ (here $\Omega_n = \{i_1, \dots, i_{m_n}\}$), then \eqref{alpha} becomes
\begin{equation}
f(x) = \sum_{\gamma \in \Gamma} {\overline{\alpha}}_{\gamma} \ggr(x+\gamma),
\end{equation}
and using \eqref{homog-fi}  we obtain
\begin{equation} \label{f-homog}
f(x) = \sum_{\gamma \in \Gamma} \beta_{\gamma}\hh(x+\gamma) \quad
\text{where} \quad \beta_{\gamma} = {\overline{\alpha}}_{\gamma} \BB_n^{-1}.
\end{equation}

We will now prove Proposition~\ref{lambda0}.
For this, let $r_0$ be the power of $x$ in $\chi_{T_n}$ (c.f.~\eqref{charact}).
Choose $m\geq n$ large enough such that
\begin{equation}\label{large-omega}
\Omega_m \supset \Omega_{n_0} - \left( D + AD + \dots + A^{r_0-1}D \right).
\end{equation}
Define the matrices $\left[(T_k)_{d}\right]_{ij}
= c_{Ai-j+d}, i,j \in \Omega_k$, for any $k \in \N$, and $d \in D$. It is shown in \cite{CHM04}, that if $x \in Q$, for any $r \geq 1$ there exists $y_r \in Q$, $\gamma_r \in \Gamma$,
such that
\begin{equation}
x  = A^{-r}(y_r + \gamma_r) \quad \text{with}\quad
\gamma_r= d_r + Ad_{r-1} + \dots + A^{r-1}d_1 \label{q-expansion}
\end{equation}
where $d_i \in D, 1 \leq i \leq r$.
Therefore, if $k \geq n_0$ and $P_k$ is as in \eqref{projection}
\begin{equation} \label{self-sim-t}
P_k\Phi(x) = (T_k)_{d_1}\dots (T_k)_{d_r}P_k\Phi(A^rx - \gamma_r) \quad x \in Q.
\end{equation}

For convenience, we will call $\Omega = \Omega_m$.
\begin{lemma}
With the previous notation, if $r \leq r_0$, then for $k\in \Omega$ and $j \in \Omega_{n_0}$,  we have
\begin{equation*}
\left[T_m^r\right]_{k (j-\gamma_r)} =  \left[(T_m)_{d_1}\dots (T_m)_{d_r}\right]_{kj},
\end{equation*}
where $\gamma_r \in D + AD + \dots + A^{r-1} D$.
\end{lemma}
\begin{remark} Note that the preceding equation does not state that both matrices are equal.
\end{remark}
\begin{proof}
We will prove the Lemma by induction on $r$. Let $\gamma_r$ be as in \eqref{q-expansion}.
\begin{itemize}
\item The case $r=1$ is trivial by the definition of $(T_m)_{d_1}$.
\item $r-1 \Longrightarrow r$
Observe first that by the choice of $\Omega$,
$$ \left[T_m\right]_{u (j-\gamma_r)} = \left[(T_m)_{d_r}\right]_{u (j-A\gamma_{r-1})}, \quad u \in \Omega, j \in \Omega_{n_0}.$$
Now
\begin{align}
\left[T_m^r\right]_{k (j-\gamma_r)}  & = \sum_{u\in \Omega} \left[T_m^{r-1}\right]_{ku}\left[T_m\right]_{u (j-\gamma_r)} \notag\\
& = \sum_{u\in \Omega} \left[T_m^{r-1}\right]_{ku}\left[(T_m)_{d_r}\right]_{u (j-A\gamma_{r-1})}
\notag\\
& = \sum_{u\in \Omega_{n_0}-\gamma_{r-1}} \left[T_m^{r-1}\right]_{ku}\left[ (T_m)_{d_r}\right]_{u (j-A\gamma_{r-1})} ,
\label{dr}
\end{align}
where the last equality follows from the $\Lambda^{\prime}$-admissibility of
$\Omega_{n_0}$.
But
\begin{equation*}
\left[(T_m)_{d_r}\right]_{u (j-A\gamma_{r-1})} = \left[(T_m)_{d_r}\right]_{(u+\gamma_{r-1}) j}, u \in \Omega_{n_0} - \gamma_{r-1}, j \in \Omega_{n_0},
\end{equation*}
and therefore, using induction and the  $\Lambda^{\prime}$-admissibility of
$\Omega_{n_0}$, \eqref{dr} becomes
\begin{align*}
\left[T_m^r\right]_{k (j-\gamma_r)}  & = \sum_{u\in \Omega_{n_0}-\gamma_{r-1}} \left[T_m^{r-1}\right]_{ku}\left[ (T_m)_{d_r}\right]_{(u+\gamma_{r-1}) j}  \\
  & = \sum_{\ell\in \Omega_{n_0}} \left[T_m^{r-1}\right]_{k (\ell -\gamma_{r-1})}\left[ (T_m)_{d_r}\right]_{\ell j} \\
  & = \sum_{\ell \in \Omega_{n_0}} \left[(T_m)_{d_1}\dots(T_m)_{d_{r-1}}\right]_{k\ell}\left[ (T_m)_{d_r}\right]_{\ell j}   \\
 & = \sum_{\ell \in \Omega} \left[(T_m)_{d_1}\dots(T_m)_{d_{r-1}}\right]_{k\ell}\left[ (T_m)_{d_r}\right]_{\ell j}
 \end{align*}
which completes the inductive step.
\end{itemize}
\end{proof}

 We can now prove Proposition~\ref{lambda0}.
\begin{proof}
Let $\Omega = \Omega_m$ be as before, and let $x \in Q \setminus (\partial Q \cup \bigcup_{i=1}^r A^{-i}\partial Q+A^{-i}D+\dots+A^{-1}D)$. Note that with this choice of $x$, $A^rx - \gamma_r \in Q^{\circ}$, and therefore,  if $u \not\in \Omega_{n_0}$, then $\varphi(A^rx - \gamma_r + u) = 0$.

Using this, together with the previous lemma, and equation \eqref{self-sim-t} for $k = m$, we have
\begin{align*}
\left[T_m^{r} P_m\Phi(A^rx)\right]_k  & =
\sum_{j \in \Omega} \left[T_m^{r}\right]_{kj} \varphi(A^rx + j) \\
& = \sum_{j \in \Omega} \left[T_m^{r}\right]_{k (j-\gamma_r)} \varphi(A^rx - \gamma_r + j) \\
& = \sum_{j \in \Omega}  \left[(T_m)_{d_1}\dots (T_m)_{d_r}\right]_{kj} \varphi(A^rx - \gamma_r + j)  \\
& = \varphi(x+k).
\end{align*}
Therefore, for $x \in Q \setminus (\partial Q \cup \bigcup_{i=1}^r A^{-i}\partial Q+A^{-i}D+\dots+A^{-1}D)$
\begin{align*}
h(x)  &=  \sum_{k \in \Omega_n} [v_i]_k \varphi(x+k) \\
& = \sum_{k \in \Omega_n} [v_i]_k
\sum_{j \in \Omega} \left[T_m^{r}\right]_{kj} \varphi(A^rx + j) \\
& =
\sum_{j \in \Omega} \left( \sum_{k \in \Omega_n} [v_i]_k \left[T_m^{r}\right]_{kj} \right)
\varphi(A^rx + j) = 0
\end{align*}
\end{proof}

\section{Accuracy and homogeneous polynomials}
\label{5}
In this section we will relate the previously obtained results, to the
 ``accuracy'' of a scaling function. We will use the
notation
of \cite{CHM98}.

\begin{definition}
The {\it accuracy} of $\varphi$ is the highest degree $\kappa$
such that all multivariate polynomials $q$  with degree $(q) < \kappa$
are in $S(\varphi)$.
\end{definition}

Let $x=(x_1,\ldots,x_d)\in\R^d$. With the standard multi-index
notation we write $x^{\alpha}=x_1^{\alpha_1} \cdots
x_d^{\alpha_d}$, where $\alpha=(\alpha_1,\ldots,\alpha_d)$ with
each $\alpha_i$ a nonnegative integer. Denote by
$|\alpha|=\alpha_1+\cdots+\alpha_d$. The number of multi-indices
$\alpha$ of degree $s$ is $d_s=
\begin{pmatrix}s+d-1\\d-1\end{pmatrix}$.

For each integer $s\geq 0$, we define the vector valued function
$X_{[s]}:\R^d\longrightarrow \C^{d_s}$ by
\begin{equation*}
X_{[s]}(x)=[x^{\alpha}]_{|\alpha|=s}, \quad x\in\R^d.
\end{equation*}
The ordering of the multi-indices $\alpha$ of degree $s$ is not
important as long as the same ordering is used throughout.

We will now look at  the behavior of $X_{[s]}(x)$ under the
multiplication
by an arbitrary $d\times d$ matrix $Z$  with scalar entries
$z_{i,j}$. If $|\alpha|=s$, then $(Zx)^{\alpha}$ is not in general
a monomial. Instead, it is a new polynomial of degree $s$, that is
still homogeneous, but possibly involves
all terms $x^{\beta}$ with $|\beta|=s$.

Let $Z_{[s]}=[z_{\alpha,\beta}^s]_{|\alpha|=s,|\beta|=s}$ be the
$d_s\times d_s$ matrix whose scalar entries $z_{\alpha,\beta}^s$
are defined by the equation
\begin{equation*}
\sum_{|\beta|=s}z_{\alpha,\beta}^sx^{\beta}=(Zx)^{\alpha}=\prod_{i=1}^d(z_{i,1}x_1+\cdots+z_{i,d}x_d)^{\alpha_i}.
\end{equation*}
The matrices $Z_{[s]}$ and their properties have been intensively
studied in
\cite{CHM98}, \cite{CHM99}. In particular, if $I_d$ denotes the
identity matrix in $\R^d$, we have that
$$(I_d)_{[s]} = I_{d_s}$$
and if $Z$ and $U$ are two
matrices,
$$ (ZU)_{[s]} = Z_{[s]} U_{[s]} \quad \text{and hence, if $Z$ is
invertible}
\quad (Z^{-1})_{[s]} = (Z_{[s]})^{-1}.$$
Dilation of $X_{[s]}(x)$ by $Z$ obeys the rule
\begin{equation*}
X_{[s]}(Zx)=Z_{[s]}X_{[s]}(x),
\end{equation*}
hence, if $A$ is the dilation matrix corresponding to the refinement
equation~\eqref{dil-eq},
$$X_{[s]}(A^{-1}x)={A^{-1}}_{[s]}X_{[s]}(x).$$
If $J_s$ is the Jordan form of ${A^{-1}}_{[s]}$, then there exists
an invertible $d_s\times d_s$ matrix $Q_s$ such that
$Q_s{A^{-1}}_{[s]}Q_s^{-1}=J_s$. So we have that
\begin{equation*}
Q_sX_{[s]}(A^{-1}x)=(Q_s{A^{-1}}_{[s]}Q_s^{-1})Q_sX_{[s]}(x).
\end{equation*}
Denote by $\widetilde{Q}_s(x)=Q_sX_{[s]}(x)$. Observe that
$\widetilde{Q}_s(x)={\left(\widetilde{Q}_s^1(x),\ldots,\widetilde{Q}_s^{d_s}(x)\right)}^t$
is a column vector of polynomials of degree $s$ that are
homogeneous. By the previous equation, we have that
\begin{equation*}
\widetilde{Q}_s(A^{-1}x)=J_s\widetilde{Q}_s(x).
\end{equation*}
Let $\beta$ be an  eigenvalue of ${A^{-1}}_{[s]}$ and $B$ the Jordan block of order $\ell$
associated to $\beta$, i.e.,
\begin{equation*}
B=\left(
\begin{array}{ccccc}
\beta & 0 & \hdots & 0 & 0 \\ 1 & \beta & \hdots & 0 & 0\\. & . &
\hdots & . & .\\ 0 & 0 & \hdots & \beta & 0 \\ 0 & 0 & \hdots & 1
& \beta
\end{array} \right)\quad \in \C^{\ell\times \ell}.
\end{equation*}
We write $\widetilde{Q}_B(x)$ for the vector that is the restriction
of $\widetilde{Q}_s(x)$ to the coordinates that correspond to the
block $B$, i.e, if $j,j+1,\ldots,j+\ell-1$ are the columns of $B$ in
$J_s$ then
$\widetilde{Q}_B(x)={\left(\widetilde{Q}_s^j(x),\ldots,
\widetilde{Q}_s^{j+\ell-1}(x)\right)}^t.$
Since $$\widetilde{Q}_s(A^{-1}x)=J_s\widetilde{Q}_s(x),$$ we have
\begin{equation}\label{ex51}
\widetilde{Q}_B(A^{-1}x)=B\widetilde{Q}_B(x).
\end{equation}
This relation will
enable us to show how, under the hypothesis of accuracy, we can
relate the Jordan form of
${A^{-1}}_{[s]}$ to the one of $T_n$. This relation
also gives a necessary condition for $\varphi$ to have
accuracy $\kappa$.
\begin{proposition}
Assume that $\varphi$ has accuracy $\kappa$ and that
$\{\varphi(\cdot-k)\}_{k\in\Gamma}$ are linearly independent. Let
$s<\kappa$. If $\beta$ is an eigenvalue of ${A^{-1}}_{[s]}$ and
$B$ is a Jordan block of ${A^{-1}}_{[s]}$ associated to $\beta$ of
order $\ell$, then $T_n$ has a Jordan block
associated to $\beta$ of order $\ell'$ with $\ell'\geq \ell$.
\end{proposition}
\begin{proof}
Consider
$\widetilde{Q}_B(x)=(\widetilde{Q}_B^1(x),\ldots,\widetilde{Q}_B^{\ell}(x))$.
It
follows from (\ref{ex51}) that
\begin{equation}\label{BB}
\begin{array}{lll}
\widetilde{Q}_B^1(A^{-1}x)&= & \beta \widetilde{Q}_B^1(x)\\
\widetilde{Q}_B^2(A^{-1}x)&=&\widetilde{Q}_B^1(x)+\beta
\widetilde{Q}_B^2(x)\\
\vdots & & \vdots \\
\widetilde{Q}_B^\ell(A^{-1}x)&=&
\widetilde{Q}_B^{\ell-1}(x)+\beta \widetilde{Q}_B^\ell(x).
\end{array}
\end{equation}
Since $\widetilde{Q}_B^i(x)\in {\mathcal S}(\varphi)$ for $1\leq i\leq
\ell$, we can write
\begin{equation*}
\widetilde{Q}_B^i(x)=Y^i\Phi(x),
\end{equation*}
for some infinite column vector $Y^i$.
From (\ref{BB}) we have for $2\leq i\leq \ell$,
\begin{equation*}
Y^i\Phi(A^{-1}x)=Y^{i-1}\Phi(x)+\beta Y^i\Phi(x),
\end{equation*}
which implies
\begin{equation*}
Y^iL\Phi(x)-\beta Y^i\Phi(x)=Y^{i-1}\Phi(x).
\end{equation*}
So, the linear independence of $\{\varphi(\cdot-k)\}_{k\in\Gamma}$
yields
\begin{equation}\label{yi}
Y^i(L-\beta I)=Y^{i-1}.
\end{equation}
Since $\widetilde{Q}_B^\ell(x)\in {\mathcal H}(A,\beta,\ell)$, by
Proposition~\ref{corr} we have that
$v=P_nY^\ell\in\kker(T_n -\beta I)^\ell$. Consider the
vectors $v_1=P_nY^\ell,v_2=(P_nY^\ell)(T_n -\beta
I),\ldots,v_\ell=(P_nY^\ell)(T_n -\beta I)^{\ell-1}$. Let us
show that $v_1,\ldots,v_\ell$ are linearly independent: Assume that
\begin{equation}\label{lin}
\sum_{i=1}^\ell\alpha_iv_i=0.
\end{equation}
Since
\begin{align*}
\left(\sum_{i=1}^\ell\alpha_iv_i\right)(T_n -\beta
I)^{\ell-1}&=\left(\sum_{i=1}^l\alpha_iv(T_n -\beta
I)^{i-1}\right)(T_n -\beta I)^{\ell-1}\\
&=(\sum_{i=1}^\ell\alpha_iv(T_n -\beta I)^{\ell+i-2}\\
&=\alpha_1v(T_n -\beta I)^{\ell-1},
\end{align*}
it follows from (\ref{lin}) that
\begin{equation*}
\alpha_1v(T_n -\beta I)^{\ell-1}=0.
\end{equation*}
Since for every $Y\in\ell(\Gamma),r\in\N$ we have that
\begin{equation*}
(P_nY)(T_n -\beta I)^r=P_n(Y(L-\beta I)^r ),
\end{equation*}
part 3 of Proposition~\ref{spect} tells us that
$v(T_n -\beta I)^{\ell-1}\neq 0$. Hence
$\alpha_1=0$. If we multiply each side of (\ref{lin}) by
$(T_n -\beta I)^{\ell-2}$ we see that $\alpha_2=0$.
Analogously $\alpha_3=\ldots=\alpha_\ell=0$ and therefore
$v_1,\ldots,v_\ell$ are linearly independent. This implies that we
have a Jordan block of $T_n$ associated to
$\beta$ of order at least $\ell$. We can repeat this procedure for
every Jordan block $B_1,\ldots,B_k$ of ${A^{-1}}_{[s]}$ associated
to $\beta$ of respective orders $l_1\geq l_2\geq\ldots \geq l_k$.
Let, for $1\leq j \leq k$
\begin{equation*}
\widetilde {Q}_{B_j}^{l_j}(x)=Y^{l_j}_j\Phi(x).
\end{equation*}
All we have to prove now is that
\begin{equation*}
\begin{array}{cccc}
P_nY^{l_1}_1,& (P_nY^{l_1}_1)(T_n -\beta I),&\ldots,&
(P_nY^{l_1}_1)(T_n -\beta I)^{l_1-1},\\
 & &\vdots& \\
P_nY^{l_k}_k,& (P_nY^{l_k}_k)(T_n -\beta I), &
\ldots,& (P_nY^{l_k}_k)(T_n -\beta I)^{l_k-1}
\end{array}
\end{equation*}
are linearly independent. Let
\begin{equation}\label{lin2}
\begin{array}{cccc}
\alpha^1_1P_nY^{l_1}_1+& \alpha^2_1(P_nY^{l_1}_1)(T_n
-\beta I)+&\ldots +&
\alpha^{l_1}_1(P_nY^{l_1}_1)(T_n -\beta
I)^{l_1-1}+\\
 & &\vdots& \\
\alpha^1_kP_nY^{l_k}_k+& \alpha^2_k(P_nY^{l_k}_k)(T_n
-\beta I)+ & \ldots+&
\alpha^{l_k}_k(P_nY^{l_k}_k)(T_n -\beta
I)^{l_k-1}=0
\end{array}
\end{equation}
Let $B_1,\ldots,B_t$ the Jordan blocks of order $l_1$. If we
multiply each side of the previous equation by
$(T_n -\beta I)^{l_1-1}$, we obtain
\begin{equation*}
\sum_{i=1}^t\alpha_i^1(P_nY^{l_1}_i)(T_n -\beta
I)^{l_1-1}=0,\text{  i.e.}
\end{equation*}
\begin{equation*}
P_n\left(\sum_{i=1}^t\alpha_i^1Y^{l_1}_i(L -\beta
I)^{l_1-1}\right)=0.
\end{equation*}
Since $\sum_{i=1}^t\alpha_i^1Y^{l_1}_i(L-\beta
I)^{l_1-1}\in\kker(L-\beta I)$, part 3 of Proposition~\ref{spect}
implies that
\begin{equation*}
\sum_{i=1}^t\alpha_i^1Y^{l_1}_i(L-\beta I)^{l_1-1}=0.
\end{equation*}
So, since by (\ref{yi}) and Proposition~\ref{li},
$Y^{l_1}_1(L-\beta I)^{l_1-1},\ldots,Y^{l_1}_t(L-\beta I)^{l_1-1}$
are linearly independent, it follows that
$\alpha_1^1=\ldots=\alpha_t^1=0$. Repeating a similar argument for
every $l_j, 2\leq j \leq k$ we can see that every scalar of
(\ref{lin2}) is equal to zero. This completes the proof.
\end{proof}
Let us now recall \eqref{ex51}, and notice that
\begin{equation*}
\widetilde{Q}_B(A^{-1}x)-\beta\widetilde{Q}_B(x)=(B-\beta
I)\widetilde{Q}_B(x).
\end{equation*}
Equivalently, if we recall the definition of $D_A$ of the previous
section,
$D_A(f)(x)=f(A^{-1}x),$ we have
\begin{equation*}
(D_A-\beta I)\widetilde{Q}_B(x)=(B-\beta I)\widetilde{Q}_B(x),
\end{equation*}
where the product on the left side is understood coordinatewise.
Moreover, for $k\in\N$,
\begin{align*}
(B-\beta I)^k\widetilde{Q}_B(x) &=\sum_{i=0}^k
\begin{pmatrix}k\\i\end{pmatrix}
 (-\beta)^{k-i}B^{i}\widetilde{Q}_B(x)\\
 &=\sum_{i=0}^k \begin{pmatrix}k\\i\end{pmatrix}
 (-\beta)^{k-i}D_A^i\widetilde{Q}_B(x)\\
 &=(D_A-\beta I)^k \widetilde{Q}_B(x).
\end{align*}
In particular, since $(B-\beta I)$ is nilpotent of order $\ell$,
we have
\begin{equation*}
(D_A-\beta I)^{\ell}\widetilde{Q}_B(x)=(B-\beta
I)^{\ell}\widetilde{Q}_B(x)=0.
\end{equation*}
Hence, all entries of $\widetilde{Q}_B(x)$ belong to ${\mathcal
H}(A,\beta,\ell)$. We can repeat this argument for every Jordan
block associated to $\beta$ and every eigenvalue $\beta$ of
${A^{-1}}_{[s]}$. It follows that each component of
$\widetilde{Q}_s(x)$ belongs to ${\mathcal H}(A,\lambda,r)$ for
some eigenvalue $\lambda$ of ${A^{-1}}_{[s]},$ and some $r\in\N$.
Since $Q_s$ is an invertible matrix and the monomials $x^{\alpha}$
with $|\alpha|=s$ are linearly independent, it follows that
$\widetilde{Q}_s^1(x),\ldots,\widetilde{Q}_s^{d_s}(x)$ are
linearly independent, and all homogeneous polynomials
$q(x)=q(x_1,\ldots,x_n)$ with $deg(q)=s$, are a linear combination
of $\widetilde{Q}_s^1(x),\ldots,\widetilde{Q}_s^{d_s}(x)$.

We can now state the next theorem:
\begin{theorem}
Assume that $\varphi$ has accuracy $\kappa$ and that
$\{\varphi(\cdot-k)\}_{k\in\Gamma}$ are linearly independent. If
$q$ is a homogeneous polynomial in $\R^d$ with $deg(q)<
\kappa$, then $q\in{\mathcal H}=\bigoplus_{\lambda \in \Delta_n}
 \mathcal{H}_{\lambda}(\varphi)$, where $\Delta_n$ is the set of
eigenvalues
 of $T_n$.
\end{theorem}
\begin{proof}
Let $s< \kappa$, and let $\widetilde{Q}_s$ and
$\widetilde{Q}_B$ be as before.
Since $\varphi$ has accuracy $\kappa$, and $s<\kappa$, all components
of
$\widetilde{Q}_B$ (in fact all components of $\widetilde{Q}_s)$
are in $S(\varphi)$, and satisfy
\begin{equation}\label{B1}
\widetilde{Q}_B(A^{-1}x)=\left(\begin{array}{ccccc} \beta & 0 & \hdots
& 0 & 0 \\ 1 & \beta & \hdots & 0 & 0\\. & . & \hdots & . & .\\ 0
& 0 & \hdots & \beta & 0 \\ 0 & 0 & \hdots & 1 & \beta
\end{array}\right)\widetilde{Q}_B(x).
\end{equation}
If we denote by $\widetilde{Q}_B^1(x)$ the first coordinate of
$\widetilde{Q}_B(x)$ we see that $\widetilde{Q}_B^1(x)$ is actually
of class
$\mathcal{H}(A,\beta,1)$.
Hence, by Proposition~\ref{corr},
$\widetilde{Q}_B^1(x)=Y\Phi$, where $P_nY\in\kker(T_n
-\beta I)$. This means that $\beta$ is also an eigenvalue of
$T_n$ and the theorem follows.
\end{proof}

The following corollary imposes conditions on the eigenvalues of
$T_n$, under the hypothesis of accuracy.
\begin{corollary}
Assume that $\varphi$ has accuracy $\kappa$ and that
$\{\varphi(\cdot-k)\}_{k\in\Gamma}$ are linearly independent. Let
$\lambda_1,\ldots,\lambda_d$ be the eigenvalues of $A$ (counted
with multiplicity). If
$\eta=(\frac{1}{\lambda_1},\ldots,\frac{1}{\lambda_d})$, then
$[\eta^{\alpha}]_{|\alpha|=s}$ are eigenvalues of
$T_n$, for $s = 0, 1, \dots, \kappa-1$.
\end{corollary}
\begin{proof}
Let  $\lambda_1,\ldots,\lambda_d$ be the eigenvalues of $A$. By
\cite{CHM98}, $[\lambda^{\alpha}]_{|\alpha|=s}$ are the
eigenvalues of $A_{[s]}$. Also, recall that since $A$ is invertible,
$A_{[s]}$ is also invertible and $(A_{[s]})^{-1}={A^{-1}}_{[s]}$. So
the eigenvalues of
${A^{-1}}_{[s]}$ are
$[\eta^{\alpha}]_{|\alpha| = s}.$ We have already
proved that if $\varphi$ has accuracy $\kappa$ and $s<\kappa$,
then every eigenvalue of ${A^{-1}}_{[s]}$ is also an eigenvalue of
$T_n$. So the result follows.
\end{proof}

\section{1-dimensional examples}

We conclude the paper by exhibiting two examples of
$(\lambda,k)$-homogeneous functions, associated to scaling
functions in dimension 1. The higher dimensional examples can be
obtained in a similar way.

\subsection{Daubechies D$_4$}
Daubechies wavelets, are those refinable functions of $N$
coefficients, that are orthogonal and provide the highest order of
accuracy possible. (Note that the splines do not form an
orthonormal base).

D$_4$ is the refinable function that satisfies the refinement
equation of $4$ coefficients:
\begin{equation} \textstyle
D_4(x) = \frac{1+\sqrt{3}}{4} D_4(2x)
+\frac{3+\sqrt{3}}{4}D_4(2x-1) + \frac{3-\sqrt{3}}{4} D_4(2x-2) +
\frac{1-\sqrt{3}}{4} D_4(2x-3).
\end{equation}
D$_4$ has accuracy 2 (it reproduces the constant and the linear
functions).

In this case the matrix $T$ has eigenvalues $1$, $\frac{1}{2}$ and
$c_0 = \frac{1+\sqrt{3}}{4}$ (see \cite{CHnM05a}). So a basis for
$\text{span}\{D_4(x), D_4(x+1), D_4(x+2)\}_{x \in [0,1]}$ is also
given by $\text{span}\{1, x, h_{c_0}(x)\}_{x \in [0,1]}$ where
$h_{c_0}$ is the homogeneous function associated to $c_0$.
\begin{figure}
\centerline {
\includegraphics[width=50mm]{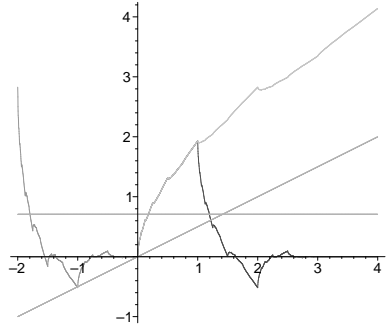}
}
\caption{\label{fig2} Daubechies D4 with the homogeneous functions. }
\end{figure}

\subsection{$(\lambda,2)$-Homogeneous function}

\begin{figure}
 \centerline{
\includegraphics[width=30mm]{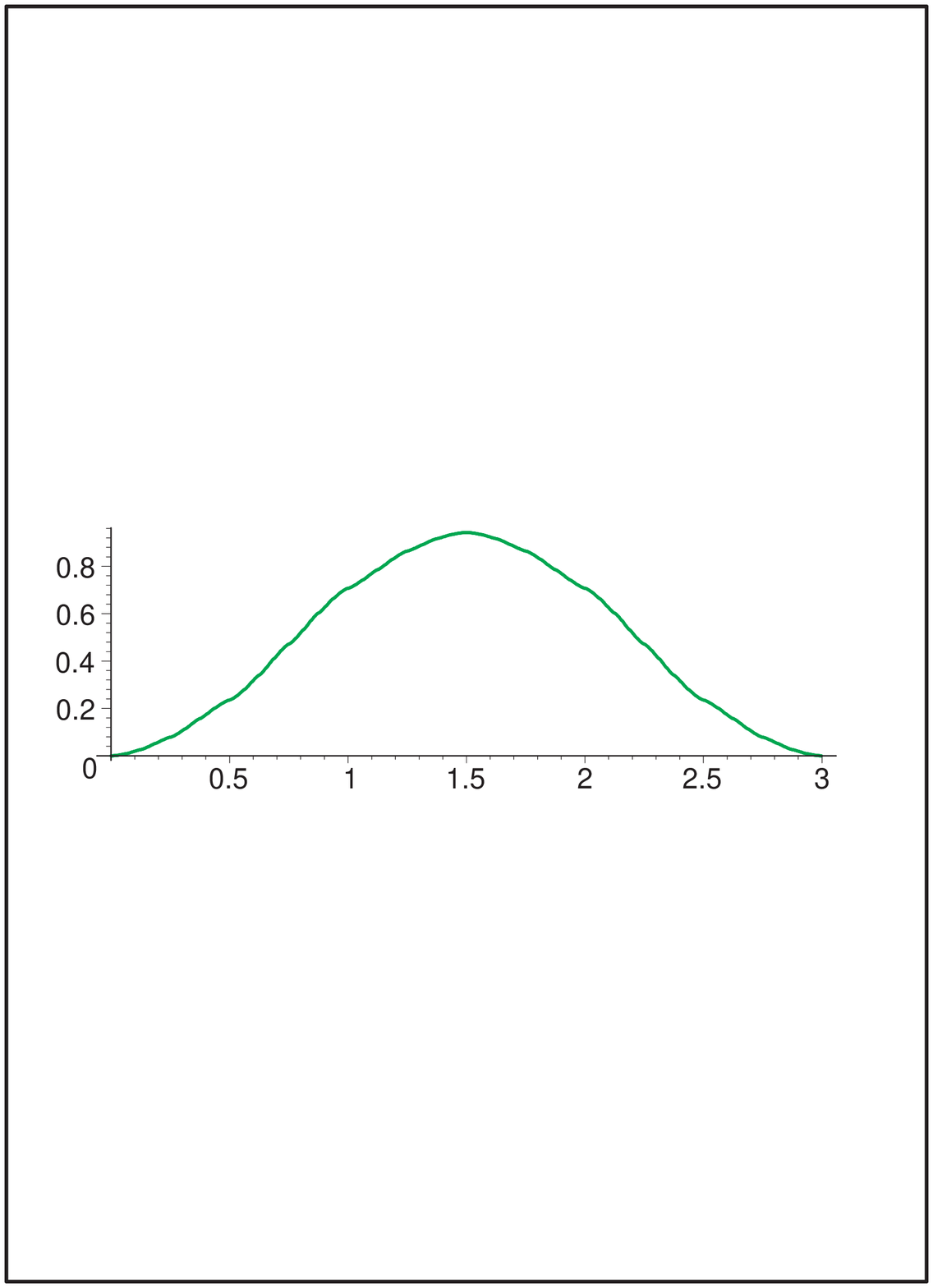}
\hspace{2cm}
\includegraphics[width=30mm]{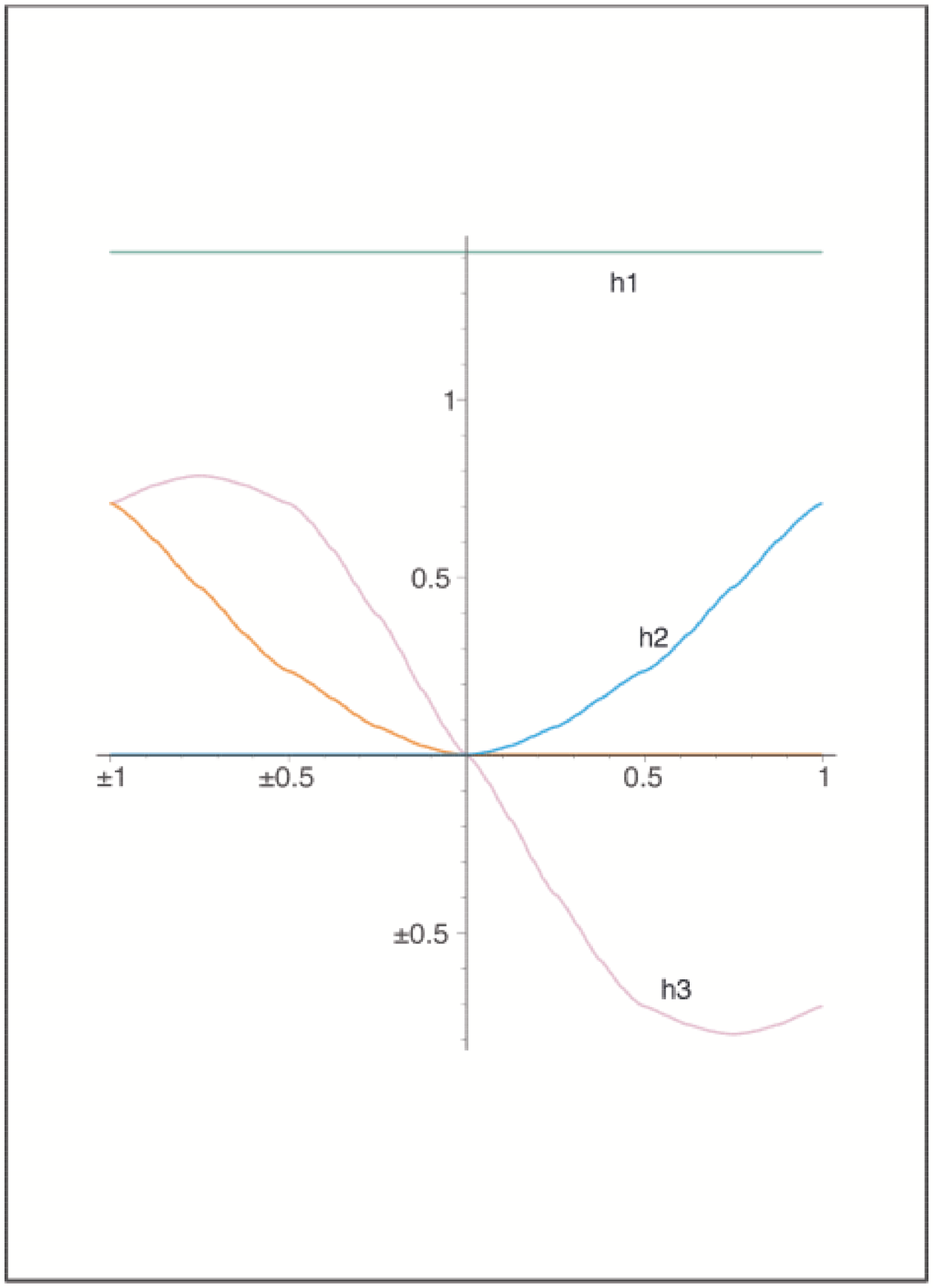}
}
\caption{\label{fig3} Scaling function for coefficients $1/3,2/3,2/3,1/3$ with Homogeneous functions - $h_1, h_2, h_3$. $h_3$ is a $2$-homogeneous function}
\end{figure}

In the previous example, we obtained a local basis of
$\text{span}\{f(x), f(x+1), f(x+2)\}_{x \in [0,1]}$ just by using
1-homogeneous functions. The following example is to illustrate,
that even in the simple case of only 4 coefficients, it may be
necessary to use homogeneous functions of order bigger than 1.
Consider the function:
\begin{equation}
f(x) = \frac{1}{3} f(2x) + \frac{2}{3}f(2x-1) + \frac{2}{3} f(2x-2)
+ \frac{1}{3} f(2x-3).
\end{equation}
It can be shown that $f$ has accuracy 1, and the eigenvalues of
$T$ are $\{1,\frac{1}{3}\}$. So in this case, $\text{span}\{f(x),
f(x+1), f(x+2)\}_{x \in [0,1]} = \text{span}\{1, h_{\{1/3,1\}}(x),
h_{\{1/3,2\}}(x)\}_{x \in [0,1]}$, where $h_{\{1/3,1\}}$ is a
1-homogeneous function corresponding to the eigenvalue $1/3$, and
$h_{\{1/3,2\}}$ is a 2-homogeneous function corresponding to the
eigenvalue $1/3$ (see Figure~\ref{fig3}).
\nopagebreak

\bibliographystyle{amsalpha}
\bibliography{cyu}

\end{document}